\documentclass[12pt]{article}
\usepackage[comma,numbers,sort&compress]{natbib}
\usepackage{rotating}
\usepackage{color}
\usepackage{amssymb}
\usepackage{amsmath,amsthm}
\usepackage{amsfonts}
\usepackage{graphicx}
\usepackage{graphics}
\usepackage{url} 
\usepackage{setspace}
\usepackage{caption}
\usepackage{lipsum}
\numberwithin{equation}{section}

\newenvironment{Proof}{\noindent{ \textbf{Proof.}}}{\hfill{$\blacksquare$}}

\usepackage{lipsum}

\begin{document}
\date{\small\textsl{\today}}
\title{Numerical Solution of HCIR Equation with Transaction Costs using Alternating Direction Implicit Method }

\author{\large
Elham Mashayekhi$^{{a}}\vspace{.5cm}$,\large
Javad Damirchi$^{{a}}$, \large 
Ahmad Reza Yazdanian$^{{b}}$\footnote{Corresponding author: yazdanian@khu.ac.ir}, \large
\\\
\small{\em $^{\mbox{\footnotesize a}}$\em Faculty of Mathematics, Statistics and Computer Science, Semnan University, Semnan, Iran}\\
\small{\em $^{\mbox{\footnotesize b}}$\em Faculty of Financial Sciences, Kharazmi University, Tehran, Iran}\\
\small{\em }\vspace{-1mm}}
\maketitle
\vspace{.1cm} 
\begin{abstract}
For valuing European options, a straightforward model is the well-known Black-Scholes formula. Contrary to market reality, this model assumed that interest rate and volatility are constant. To modify the Black-Scholes model, Heston and Cox-Ingersoll-Ross (CIR) offered the stochastic volatility and the stochastic interest rate models, respectively. The combination of the Heston, and the CIR models is called the Heston-Cox-Ingersoll-Ross (HCIR) model. Another essential issue that arises when purchasing or selling a good or service is the consideration of transaction costs which was ignored in the Black-Scholes technique. Leland improved the simplistic Black-Scholes strategy to take transaction costs into account. The main purpose of this paper is to apply the alternating direction implicit (ADI) method at a uniform grid for solving the HCIR model with transaction costs in the European style and comparing it with the explicit finite difference (EFD) scheme. Also, as evidence for numerical convergence, we convert the HCIR model with transaction costs to a linear PDE (HCIR) by ignoring transaction costs, then we estimate the solution of HCIR PDE using the ADI method which is a class of finite difference schemes, and compare it with analytical solution and EFD scheme. For multi-dimensional Black-Scholes equations, the ADI method, which is a category of finite difference techniques, is appropriate. When the dimensionality of the space increases, finite difference techniques frequently become more complex to perform, comprehend, and apply. Consequently, we employ the ADI approach to divide a multi-dimensional problem into several simpler, quite manageable sub-problems to overcome the dimensionality curse.\\

\vspace{.5cm}
\end{abstract}
\begin{quotation}
\section{Introduction}
\label{section1}
The Black-Scholes method \cite{B_1973}, in which the underlying asset exhibits geometric Brownian motion, is the standard model used to calculate the price of the European option. The assumption that the interest rate and volatility are constant in this model is contrary to market reality. To improve results, a modified Black-Scholes strategy including stochastic volatility was advocated by many authors (see \cite{H_W_1987,S_Z_1999,S_S_1991}). To work around a Black-Scholes model restriction, Heston \cite{Heston_1993} developed the stochastic volatility approach in 1993.\\
The construction of interest models is essential because, in addition to volatility, interest rates also have an impact on the values of many practically significant assets, such as callable or path-dependent securities. Some research \cite{G_O_W_2011,G_O_W_2012,G-G_O_2013,H_L_P_S_2000,C_I_R_1985,H_W_1990,V_1977} described the interest rate and volatility as stochastic. A modified Black-Scholes problem with stochastic interest rate and stochastic volatility models, the Heston-Cox-Ingersoll-Ross (HCIR) model was introduced by Sippel and Ohkoshi \cite{S_O_2002}. These problems do not have a closed-form solution, so researchers used numerical methods to estimate the solution. For option pricing in the general Black-Scholes model, several numerical methods were proposed, including the Fourier techniques \cite{G_O_2011}, the finite difference (FD) methods \cite{W_2019,C_W_Z_2021,D_2006}, the splitting methods \cite{H_2013}, the binomial method \cite{Z_2002}, the Monte Carlo method \cite{N_Y_2011}, and so on. 
For multi-dimensional Black-Scholes equations, the ADI method, which is a category of finite difference techniques, is appropriate. When the dimensionality of the space increases, finite difference techniques frequently become more complex to perform, comprehend, and apply. Consequently, we employ the ADI approach to divide a multi-dimensional issue into several simpler, quite manageable sub-problems to overcome the dimensionality curse. There are several benefits to the ADI method, which was pioneered in the United States by Douglas \cite{Douglas,D_1955}, Peaceman and Rachford \cite{Peaceman,P_R_1955}, Douglas and Rachford \cite{D_R_1956} and others. Express contrast strategies are seldom used to take care of beginning limit esteem issues inferable from their unfortunate soundness issues. Unfortunately, implicit approaches are challenging to solve in more than two dimensions despite their improved stability. Because they are programmable by resolving a simple tridiagonal system of equations, ADI methods became an alternative. During the period that the ADI method was developing Yanenko \cite{Y_1971} was developing splitting methods for resolving three-dimensional time subordinate partial differential equations. The financial literature favors the ADI method. However, it can be applied in a variety of ways to break the Black-Scholes model as one-dimensional, smaller issues. An alternating direction strategy was introduced by Douglas \cite{D_1962} to solve nonlinear parabolic and elliptic problems involving boundary conditions in a rectangular area. 
According to the findings, if the partial differential equation is solved smoothly enough; The error of solution for a linear partial differential equation is $O({(\Delta x)^2} + {(\Delta t)^2})$ , whereas the accuracy reduces to $O({(\Delta x)^2} + {\Delta t})$ when a partial differential equation with nonlinear terms is solved using the predictor-corrector generalization. Shidfar et al. \cite{S_2014} used the Peaceman-Rachford method to obtain European Spread options on two assets. They considered the first-order feedback approach, which creates the linear partial differential equation (PDE) and explains the convergence and stability.
Yazdanian and Pirvu \cite{Y_2014} also utilized this strategy to numerically solve the Spread option using the full-feedback model. The unconditionally stable high-order compact FD method used by D$\ddot u$ring and Fournie \cite{D_F_2012} to determine the option value in stochastic volatility models. This scheme is accurate to the fourth order in space and the second order in time. later, Safaei et al. \cite{S_N_N_2019} for the HCIR option pricing model used the generalized componentwise splitting strategy. In three-dimensional space, they represent the fair option price for bonds and derivatives. In the past, the multi-dimensional Black-Scholes model was numerically solved (see \cite{J_S_M_W_2022}). They used the operator splitting method (OSM) and the explicit Crank-Nicolson scheme.\\ 
The consideration of transaction costs is another crucial issue that arises when trading a commodity. The lack of any evidence to back up the assumption that transaction costs exist in actual markets is another restriction of the Black-Scholes approach. Leland \cite{Leland_1985} proposed the method for hedging a portfolio that rebalanced itself at each time step and developed the option pricing model with transaction costs in 1985.  A hedge contingent claim was explained by Hodges and Neuberger \cite{H_N_1989}, based on similar transaction costs. They were able to accurately replicate a claim, which was easier than using the Leland method. Grannan and Swindle \cite{G_S_1996} described a novel strategy with transaction costs including the Leland method. Also, transaction costs were utilized by Zhao and Ziemba \cite{Z_Z_2007,Zh_Z_2007} to estimate the Leland model-based option value and simulate volatility. An option value was determined using a stochastic interest rate and transaction costs proposed in \cite{Se_2014}. The exchange costs produced nonlinear parts in such a halfway differential condition by utilizing PDE2D programming to tackle a difficult partial differential equation that included irregular unpredictability and exchange costs. In response to applications in finance, Mariani et al. \cite{M_2015} employed software to solve a challenging partial differential equation. The solution provides an expense of the European option that takes transaction fees and unpredictable volatility into consideration. 
 In their initial work, Nguyen and Pergamenshchikov \cite{N_P_2015} included jumps and introduced the replication method for approximating option pricing under proportional transaction costs and stochastic volatility. They also demonstrated several limitation statements regarding the normalized replication error of the Leland method and established the characteristics of underhedging, including \cite{N_P_2017}.\\
Although Bakshi et al. \cite{G_B_1997} discovered an analytical solution to the HCIR problem, the HCIR model with transaction costs typically lacks a closed-form solution. In this paper, the Douglas scheme is used to estimate the European option price using the HCIR model with transaction costs. In Section \ref{section2}, we present the HCIR model with this kind of partial correlation formation and provide an analytical solution to the PDE with zero-coupon bonds. We get HCIR with transaction costs PDE in Section \ref{section3}. The Douglas approach is used to evaluate the HCIR formula for estimating the price of the European option in Section \ref{section4}. The numerical results of the HCIR with transaction costs PDE is presented in Section \ref{section5}. We estimate the solution of the HCIR PDE using the Douglas method and compare it with the analytical solution and explicit finite difference method to demonstrate numerical stability. We do this by converting the nonlinear PDE for the HCIR with transaction costs to the linear PDE and omitting the transaction costs. The remarks and conclusion are then discussed. 
\section{The HCIR Model and Pricing Zero-Coupon Bonds}
\label{section2}
We consider the HCIR approach with the risk-neutral measure $Q$ as displayed in (see \cite{C_J_L_2016})
\begin{equation}
dS_t = R_tS_tdt + \sqrt {V_t} S_tdW_t^1,
\label{(S1_2)}
\end{equation}
\begin{equation}
dV_t = k\left( {\zeta  - V_t} \right)dt + \sigma \sqrt {V_t} dW_t^2,
\label{(S2_2)}
\end{equation}
\begin{equation}
dR_t = a \left( {b  - R_t} \right)dt + \eta \sqrt {R_t} dW_t^3.
\label{(S3_2)}
\end{equation}
The parameters are the positive constants $k$, $\zeta$, $\sigma$, $\eta$, $a$, and $b$, and $S_t$, $V_t$, and $R_t$ correspond to the asset values, related volatility, and such interest rate at time $t>0$, respectively. The standard Brownian movements are depicted by $\{ W_t^1:t \geqslant 0\}, \{ W_t^2:t \geqslant 0\}$ and $\{ W_t^3:t \geqslant 0\}$ while utilizing risk-neutral measure $Q$. A Heston stochastic volatility model \cite{Heston_1993} is displayed in Eq. (\ref{(S2_2)}). The Cox-Ingersoll-Ross (CIR) model Eq. (\ref{(S3_2)}) has always been used to describe interest rates because it uses a physical process.
In Eqs. (\ref{(S1_2)})-(\ref{(S3_2)}), the following relationships are presumed (see \cite{W_2019}): 
\begin{equation}
\label{(S3_22)}
 < dW_t^1,dW_t^2 >  = \rho dt,\begin{array}{*{20}{c}}
  {}&{ < dW_t^1,dW_t^3 >  = 0,}&{ < dW_t^2,dW_t^3 >  = 0.}
\end{array}
\end{equation}
\\
We assume that the discount process is as follows (see \cite{W_2019}):
\begin{equation}
D(t) ={e^{ - \int_0^t {R(\tau )d\tau } }}
\end{equation}
Consequently, the method used to calculate the cost of accounts for money markets is
\begin{equation}
\begin{gathered}
  \frac{1}{{D(t)}} = {e^{\int_0^t {R(\tau )d\tau } }} \hfill \\
  dD(t) =  - R(t){e^{ - \int_0^t {R(\tau )d\tau } }}dt =  - R(t)D(t)dt. \hfill \\ 
\end{gathered}
\end{equation}
A bond with a zero coupon and a maturity at time $T$ is valued $Y(t,T)$ as a time of $t$. Since the discounted price of such a bond needs to be a martingale based on the risk-neutral measurement, let's suggest that $Y(T,T)=1$. 
\begin{equation}
{\rm E}_t^Q(Y(T,T)D(T)) = {\rm E}_t^Q(D(T)) = D(t)Y(t,T).
\end{equation}
By using the equation as a guide, we get  
\begin{equation}
Y(t,T) = {\rm E}_t^Q\left( {{e^{ - \int_t^T {R(\tau )d\tau } }}} \right)
\end{equation}
Because $dR$ should be a Markov chain with  
\begin{equation}
Y(t,T) = Z(t,R(t)).
\end{equation}
The PDE for $Z(t,R(t))$ is obtained by differentiating $D(t)Z(t,R(t))$ using It$\widehat o$'s lemma, 
\begin{equation}
\begin{gathered}
  d\left( {D(t)Z(t,R(t))} \right) = Z(t,R(t))dD(t) + D(t)dZ(t,R(t)) \hfill \\
   = D(t)\left( { - RZdt + \frac{{\partial Z}}{{\partial t}}dt + \frac{{\partial Z}}{{\partial R}}dR + \frac{1}{2}\frac{{{\partial ^2}Z}}{{\partial {R^2}}}d{R^2}} \right) \hfill \\ 
\end{gathered}
\end{equation}
After simplification, we have (following \cite{W_2019})
\begin{equation}
\frac{{\partial Z}}{{\partial t}} + \frac{1}{2}{\eta ^2}R\frac{{{\partial ^2}Z}}{{\partial {R^2}}} + a((b - R) - \lambda \eta )\frac{{\partial Z}}{{\partial R}} - RZ = 0,\begin{array}{*{20}{c}}
  {}&{Z(R,T;T) = 1,} 
\end{array}
\label{(S_2)}
\end{equation}
where $\lambda \sqrt R$ corresponds to the market price of risk.
According to Shreve and Steven \cite{SH_2004}, the closed form solution for PDE (\ref{(S_2)}) is
\begin{equation}
\ln \left[ {Z(R,t;T)} \right] =  - {B_1}R + {B_2},
\label{(S5_2)}
\end{equation}
Such that 
\begin{equation}
{B_1} =\frac{{2\alpha \beta }}{{{\eta ^2}}}\ln \left( {\frac{{2\gamma {e^{(\gamma  + \alpha )(T - t)/2}}}}{{(\alpha  + \gamma )({e^{\gamma (T - t)}} - 1) + 2\gamma }}} \right),
\label{(S6_2)}
\end{equation}
\begin{equation}
{B_2} = \frac{{2\left( {{e^{\gamma (T - t)}} - 1} \right)}}{{(\alpha + \gamma )({e^{\gamma (T - t)}} - 1) + 2\gamma }},
\label{(S7_2)}
\end{equation}
\begin{equation}
\gamma  = \sqrt {{\alpha^2} + 2{\eta ^2}.} 
\label{(S8_2)}
\end{equation}
\section{The HCIR Model Including Transaction Costs}
\label{section3}
We took into account the quantity $C(S,V,R,t)$ additionally the amounts $-{\Delta _1} , - {\Delta _2}$ and $-{\Delta _3}$ of such asset price, variance, and zero-coupon bond respectively that taken into consideration as a portfolio of $\Pi$. The effects of using a zero-coupon bond to hedge the stochastic interest rate are as follows (following \cite{W_2019})
\begin{equation}
\label{(S1_3)}
\Pi  = C - {\Delta _1}S - {\Delta _2}V - {\Delta _3}Z.
\end{equation}
The self-financing argument states that after taking into account transaction costs and stochastic interest rates, we have 
\begin{equation}
\label{(S3.2)}
d\Pi  = dC - {\Delta _1}dS - {\Delta _2}dV - {\Delta _3}dZ.
\end{equation}
To determine the dynamics of $C$, we use It${\hat o}$'s formula,\\ \\
\begin{equation}
\label{(S3.3)}
\begin{gathered}
  dC = \frac{{\partial C}}{{\partial t}}dt + \frac{{\partial C}}{{\partial S}}dS + \frac{{\partial C}}{{\partial R}}dR + \frac{{\partial C}}{{\partial V}}dV + \frac{1}{2}V{S^2}\frac{{{\partial ^2}C}}{{\partial {S^2}}}dt \hfill \\
   + \rho \sigma VS\frac{{{\partial ^2}C}}{{\partial S\partial V}}dt + \frac{1}{2}{\sigma ^2}V\frac{{{\partial ^2}C}}{{\partial {V^2}}}dt + \frac{1}{2}{\eta ^2}R\frac{{{\partial ^2}C}}{{\partial {R^2}}}dt \hfill \\ 
\end{gathered}
\end{equation}
\\ and
\begin{equation}
\label{(S3.4)}
  dZ = \frac{{\partial Z}}{{\partial t}}dt + \frac{{\partial Z}}{{\partial R}}dR+ \frac{1}{2}{\eta ^2}R\frac{{{\partial ^2}Z}}{{\partial {R^2}}}dt.
\end{equation}
Substituting Eqs. (\ref{(S3.3)}),(\ref{(S3.4)}) in Eq. (\ref{(S3.2)}) as follows:
\begin{equation}
\label{(S3.5)}
\begin{gathered}
  d\Pi  = \left( {\frac{{\partial C}}{{\partial t}} + \frac{1}{2}V{S^2}\frac{{{\partial ^2}C}}{{\partial {S^2}}} + \rho \sigma VS\frac{{{\partial ^2}C}}{{\partial S\partial V}} + \frac{1}{2}{\sigma ^2}V\frac{{{\partial ^2}C}}{{\partial {V^2}}} + \frac{1}{2}{\eta ^2}R\frac{{{\partial ^2}C}}{{\partial {R^2}}} - {\Delta _2}\frac{{\partial Z}}{{\partial t}} - \frac{{{\Delta _2}}}{2}{\eta ^2}R\frac{{{\partial ^2}Z}}{{\partial {R^2}}}} \right)dt \hfill \\
   + \left( {\frac{{\partial C}}{{\partial S}} - \Delta } \right)dS + \left( {\frac{{\partial C}}{{\partial V}} - {\Delta _1}} \right)dV + \left( {\frac{{\partial C}}{{\partial R}} - {\Delta _2}\frac{{\partial Z}}{{\partial R}}} \right)dR - {k_0}S{\nu _0} - {k_1}V{\nu _1} - {k_2}Z{\nu _2}, \hfill \\ 
\end{gathered}
\end{equation}
where ${k_0}S{\nu _0},{k_1}V{\nu _1},{k_2}Z{\nu _2}$ represent the transaction costs of trading the asset price, the volatility, and the zero-coupon bond, respectively. Eq. (\ref{(S3.5)}) fulfills the partial differential equation (see Wang \cite{W_2019}):
\begin{equation}
\label{(S2_18)}
\begin{gathered}
  \frac{{\partial C}}{{\partial t}} =  - \frac{1}{2}V{S^2}\frac{{{\partial ^2}C}}{{\partial {S^2}}} - \frac{1}{2}V{\sigma ^2}\frac{{{\partial ^2}C}}{{\partial {V^2}}} - \rho \sigma VS\frac{{{\partial ^2}C}}{{\partial S\partial V}} - \frac{1}{2}{\eta ^2}R\frac{{{\partial ^2}C}}{{\partial {R^2}}} - RS\frac{{\partial C}}{{\partial S}} - RV\frac{{\partial C}}{{\partial V}} \hfill \\
   - \alpha (\beta  - R)\frac{{\partial C}}{{\partial R}} + RC + {F_1} + {F_2} + {F_3}, \hfill \\ 
\end{gathered}
 \end{equation}\\
where\\
$\begin{gathered}
  {F_1} = {k_0}S\sqrt {\frac{{2}}{\pi \delta t }} \sqrt {V{S^2}{{\left( {\frac{{{\partial ^2}C}}{{\partial {S^2}}}} \right)}^2} + {\sigma ^2}V{{\left( {\frac{{{\partial ^2}C}}{{\partial S\partial V}}} \right)}^2} + 2\rho V\sigma S\left( {\frac{{{\partial ^2}C}}{{\partial {S^2}}}} \right)\left( {\frac{{{\partial ^2}C}}{{\partial S\partial V}}} \right)} , \hfill \\
  {F_2} = {k_1}V\sqrt {\frac{{2}}{\pi \delta t }} \sqrt {{\sigma ^2}V{{\left( {\frac{{{\partial ^2}C}}{{\partial {V^2}}}} \right)}^2} + V{S^2}{{\left( {\frac{{{\partial ^2}C}}{{\partial S\partial V}}} \right)}^2} + 2\rho V\sigma S\left( {\frac{{{\partial ^2}C}}{{\partial {V^2}}}} \right)\left( {\frac{{{\partial ^2}C}}{{\partial S\partial V}}} \right)} , \hfill \\
  {F_3} = {k_2}\sqrt R \sqrt {\frac{{2}}{\pi \delta t }} \frac{Z}{{\left| \theta  \right|}}\eta \left| {\frac{{{\partial ^2}C}}{{\partial {R^2}}}} \right|,\theta  = \frac{{\partial Z}}{{\partial R}}
\end{gathered}$ \\
and boundary circumstances such 
\begin{equation}
\begin{gathered}
  C(0,V,R,t) = 0,\frac{{\partial C}}{{\partial S}}({S_{\max }},V,R,t) = 1, \hfill \\
  \frac{{\partial C}}{{\partial t}}(S,0,R,t) =-\frac{1}{2}{\eta ^2}R\frac{{{\partial ^2}C}}{{\partial {R^2}}}-RS\frac{{\partial C}}{{\partial S}}(S,0,R,t) - \alpha(\beta - R)\frac{{\partial C}}{{\partial R}}(S,0,R,t) + RC(S,0,R,t)  \hfill \\
+ {F_3}(S,0,R,t) , \hfill \\
  C(S,{V_{\max }},R,t) = S, \hfill \\
  \frac{{\partial C}}{{\partial R}}(S,V,0,t) = 0,\frac{{\partial C}}{{\partial R}}(S,V,{R_{\max }},t) = 0. \hfill \\ 
\end{gathered}
 \end{equation}\\
We define the final condition for the European call option with a strike price of $E$ and a time maturity of $T$ as follows
\begin{equation}
\label{(S3_19)}
C(S,V,R,T) = \max \left( {S - E,0} \right).
 \end{equation}\\
\section{Alternating Direction Implicit method}
\label{section4}
In this part, we approximate the PDE (\ref{(S2_18)}) using the alternating direction implicit method by changing the variable $t = T - \tau$. It consists of variables $\tau \in \left[ {0,T} \right]$, $R \in \left[ {0,{R_{\max }}} \right]$, $V \in \left[ {0,{V_{\max }}} \right]$, and $S \in \left[ {0,{S_{\max }}} \right]$. How about we expect that there are focuses $(i\Delta S,j\Delta V,k\Delta R,n\Delta \tau)$ in a uniform with the accompanying files $i=0,1,...,M$ and $j=0,1,...,J$, $k=0,1,...K$, $n=0,1,...,N$ where $M\Delta S = {S_{\max }},J\Delta V = {V_{\max }},K\Delta R = {R_{\max }}$ and  $N\Delta \tau=T$. We take $C_{i,j,k}^n$ to be the numerical solution of PDE (\ref{(S2_18)}) at point $(i\Delta S,j\Delta V,k\Delta R,n\Delta \tau)$ and discretize PDE (\ref{(S2_18)}) utilizing the accompanying finite difference scheme:\\ \\
$\begin{gathered}
  \frac{{\partial C}}{{\partial t}}({S_i},{V_j},{R_k},{t_n}) \approx \frac{{C({S_i},{V_j},{R_k},{t_n} + \Delta t) - C({S_i},{V_j},{R_k},{t_n})}}{{\Delta t}} + O(\Delta t), \hfill \\
  \frac{{\partial C}}{{\partial S}}({S_i},{V_j},{R_k},{t_n}) \approx \frac{{C({S_i} + \Delta S,{V_j},{R_k},{t_n}) - C({S_i} - \Delta S,{V_j},{R_k},{t_n})}}{{2\Delta S}} + O({(\Delta S)^2}), \hfill \\
  \frac{{\partial C}}{{\partial V}}({S_i},{V_j},{R_k},{t_n}) \approx \frac{{C({S_i},{V_j} + \Delta V,{R_k},{t_n}) - C({S_i},{V_j} - \Delta V,{R_k},{t_n})}}{{2\Delta V}} + O({(\Delta V)^2}), \hfill \\
  \frac{{\partial C}}{{\partial R}}({S_i},{V_j},{R_k},{t_n}) \approx \frac{{C({S_i},{V_j},{R_k} + \Delta R,{t_n}) - C_{i,j,k - 1}^n({S_i},{V_j},{R_k} - \Delta R,{t_n})}}{{2\Delta R}} + O({(\Delta R)^2}), \hfill \\
 \frac{{{\partial ^2}C}}{{\partial {S^2}}}({S_i},{V_j},{R_k},{t_n}) \approx \frac{{C({S_i} + \Delta S,{V_j},{R_k},{t_n}) - 2C({S_i},{V_j},{R_k},{t_n}) + C({S_i} - \Delta S,{V_j},{R_k},{t_n})}}{{{{(\Delta S)}^2}}}  \hfill \\
+ O({(\Delta S)^2}), \hfill \\
  \frac{{{\partial ^2}C}}{{\partial {V^2}}}({S_i},{V_j},{R_k},{t_n}) \approx \frac{{C({S_i},{V_j} + \Delta V,{R_k},{t_n}) - 2C({S_i},{V_j},{R_k},{t_n}) + C({S_i},{V_j} - \Delta V,{R_k},{t_n})}}{{{{(\Delta V)}^2}}}  \hfill \\
+ O({(\Delta V)^2}), \hfill \\
  \frac{{{\partial ^2}C}}{{\partial {R^2}}}({S_i},{V_j},{R_k},{t_n}) \approx \frac{{C({S_i},{V_j},{R_k} + \Delta R,{t_n}) - 2C({S_i},{V_j},{R_k},{t_n}) + C({S_i},{V_j},{R_k} - \Delta R,{t_n})}}{{{{(\Delta R)}^2}}}  \hfill \\
+ O({(\Delta R)^2}), \hfill \\
  \begin{gathered}
  \frac{{\partial C}}{{\partial S\partial V}}({S_i},{V_j},{R_k},{t_n}) \approx \frac{1}{{4\Delta S\Delta V}}(C({S_i} + \Delta S,{V_j} + \Delta V,{R_k},{t_n}) + C({S_i} - \Delta S,{V_j} - \Delta V,{R_k},{t_n})   \hfill \\
   - C({S_i} + \Delta S,{V_j} - \Delta V,{R_k},{t_n}) - C({S_i} - \Delta S,{V_j} + \Delta V,{R_k},{t_n})) + O({(\Delta S)^2} + {(\Delta V)^2}). \hfill \\ 
\end{gathered} \hfill \\  \\
\end{gathered}$
We consider the following operator
\begin{equation}
\label{(Star_2)}
A={A_{SV}} + {A_S} + {A_V} + {A_R},
\end{equation}
where
\begin{equation}
\label{(SS4_2)}
\begin{gathered}
  {A_S}C = \frac{1}{2}V{S^2}\frac{{{\partial ^2}C}}{{\partial {S^2}}} + RS\frac{{\partial C}}{{\partial S}} - \frac{R}{3}C, \hfill \\
  {A_V}C = \frac{1}{2}V{\sigma ^2}\frac{{{\partial ^2}C}}{{\partial {V^2}}} + RV\frac{{\partial C}}{{\partial V}} - \frac{R}{3}C, \hfill \\
  {A_R}C = \frac{1}{2}{\eta ^2}r\frac{{{\partial ^2}C}}{{\partial {R^2}}} + \alpha(\beta - R)\frac{{\partial C}}{{\partial R}} - \frac{R}{3}C, \hfill \\
  {A_{SV}}C = \rho \sigma VS\frac{{{\partial ^2}C}}{{\partial S\partial V}}. \hfill \\ 
\end{gathered}
\end{equation}
By numerically solving PDE (\ref{(S2_18)}) using the Douglas method, as follows (see \cite{H_2013,D_1962}):
\begin{equation}
\label{(S4_2)}
\begin{gathered}
  \frac{{{C({S_i},{V_j},{R_k},{t_n})} - {C({S_i},{V_j},{R_k},{t_{n-1}})}}}{{\Delta \tau}} = \theta_1 {A_S}{C({S_i},{V_j},{R_k},{t_n})} + \left( {1 - \theta_1 } \right){A_S}{C({S_i},{V_j},{R_k},{t_{n-1}})} \hfill \\ 
+ \theta_1 {A_V}{C({S_i},{V_j},{R_k},{t_n})} + (1 - \theta_1 ){A_V}{C^C({S_i},{V_j},{R_k},{t_{n-1}})}+ \theta_1 {A_R}{C({S_i},{V_j},{R_k},{t_n})}  \hfill \\
 + \left( {1 - \theta_1 } \right){A_R}{C({S_i},{V_j},{R_k},{t_{n-1}})} + {A_{SV}}{C({S_i},{V_j},{R_k},{t_{n-1}})} + {\varphi ({S_i},{V_j},{R_k},{t_{n-1}})}, \hfill \\ 
\end{gathered}
\end{equation}
that ${\varphi ({S_i},{V_j},{R_k},{t_{n-1}})}$ includes nonlinear terms of PDE  (\ref{(S2_18)}) at point $(i\Delta S,j\Delta V,k\Delta R,(n - 1)\Delta \tau)$.
Multiply Eq. (\ref{(S4_2)}) by ${{\Delta \tau}}$ and rearrange: \\
\begin{equation}
\label{(S4_3)}
\begin{gathered}
  \left( {I - \theta_1 \Delta \tau{A_S} - \theta_1 \Delta \tau{A_V} - \theta_1 \Delta \tau{A_R}} \right){C({S_i},{V_j},{R_k},{t_n})} = (I + (1 - \theta_1 )\Delta \tau{A_S} + (1 - \theta_1 )\Delta \tau{A_V} \hfill \\
 + (1 - \theta_1 )\Delta \tau{A_R}+\Delta \tau{A_{SV}}){C({S_i},{V_j},{R_k},{t_{n-1}})} + \Delta \tau{\varphi ({S_i},{V_j},{R_k},{t_{n-1}})}, \hfill \\ 
\end{gathered}
\end{equation}
where $I$ presents the identity operator. By adding $[{\theta_1 ^2}{(\Delta \tau)^2}({A_R}{A_V} + {A_V}{A_S} + {A_R}{A_S}) \hfill \\
- {\theta_1 ^3}{(\Delta \tau)^3}{A_S}{A_V}{A_R}]{C({S_i},{V_j},{R_k},{t_n})}$ on the left side and $[{\theta_1 ^2}{(\Delta \tau)^2}({A_R}{A_V} + {A_V}{A_S} + {A_R}{A_S}) - {\theta_1 ^3}{(\Delta \tau)^3}{A_S}{A_V}{A_R}]{C({S_i},{V_j},{R_k},{t_{n-1}})}$ on the right side, we have
\begin{equation}
\label{(S4_4)}
\begin{gathered}
  (I - \theta_1 \Delta \tau{A_R})(I - \theta_1 \Delta \tau{A_V})(I - \theta_1 \Delta \tau{A_S}){C({S_i},{V_j},{R_k},{t_n})} = (I - \theta_1 \Delta \tau{A_R}) \hfill \\ 
(I - \theta_1 \Delta \tau{A_V})(I - \theta_1 \Delta \tau{A_S}){C({S_i},{V_j},{R_k},{t_{n-1}})}+ \Delta \tau A{C({S_i},{V_j},{R_k},{t_{n-1}})} \hfill \\
+\Delta \tau{\varphi ({S_i},{V_j},{R_k},{t_{n-1}})}, \hfill \\ 
\end{gathered}
\end{equation}
and this gives Douglas method \cite{H_2013}
\begin{equation}
\label{(S4_1)}
\begin{gathered}
  {G_0} = (1+ \Delta \tau)A{C({S_i},{V_j},{R_k},{t_{n-1}})} +{\varphi ({S_i},{V_j},{R_k},{t_{n-1}})}, \hfill \\
  {G_1} = {G_0} + \theta_1\Delta t{A_S}\left( {{G_1} - {C({S_i},{V_j},{R_k},{t_{n-1}})}} \right), \hfill \\
  {G_2} = {G_1} + \theta_1\Delta t{A_V}\left( {{G_2} - {C({S_i},{V_j},{R_k},{t_{n-1}})}} \right), \hfill \\
  {C({S_i},{V_j},{R_k},{t_{n}})} = {G_2} + \theta_1\Delta t{A_R}\left( {{C({S_i},{V_j},{R_k},{t_{n}})} - {C({S_i},{V_j},{R_k},{t_{n-1}})}} \right), \hfill \\ 
\end{gathered}
\end{equation} \\
The unconditional stability by $\theta_1  = \frac{2}{3}$ was obtained by applying the Douglas method to convection-diffusion problems involving mixed derivative factors in the Von Neumann technique (following \cite{H_K_2013}). The following section provides the value of a European call option at $\tau=T$. 
\section{Numerical Results}
\label{section5}
In this section, we obtain numerical solutions to PDE (\ref{(S2_18)}) using the ADI method in MATLAB programming language. To estimate the solution of the PDE (\ref{(S2_18)}), the following assumptions are made: The maximum stock price is considered $5E$ with the strike price set at $100$. The maximum volatility and maximum interest rate are set as 1. Additionally, the final time is considered as 1. To approximate the solution of HCIR PDE with transaction costs, we take into account specific parameters listed in Table \ref{Table 1}. Tables \ref{Table 2} and \ref{Table 22} display the price of a European call option with $k_0=k_1=k_2=0,V=0.2$  at $R=0.2$ and $R=0.4$, respectively. Furthermore, we contrast the results with the explicit finite difference (FD) method \cite{W_2019} and also the actual solution \cite{G_B_1997}. In Figures \ref{figure1_1} and \ref{fig2_1}, we contrast an analytical solution with an approximation solution at $R=0.2,V=0.2$ and $R=0.4,V=0.4$ respectively. In Figures \ref{figure1} and \ref{fig1}, an approximation solution was created using the ADI method at $R=0.2,R=0.4$ respectively. The European call option price is written in Table \ref{Table 5} at $R=0.2,R=0.4$ with $M=1200,J=80,K=80,N=30$ , and $k_0=k_1=k_2=0.02$. The results are also compared with those obtained using the explicit finite difference method.\\
 In Figure \ref{figure2}, payoff values are lower than the numerical solution of the HCIR PDE (\ref{(S2_18)}) when compared to the European call option price with ${k_0}={k_1}={k_2}=0.02$ and $V=0.2,V=0.7$.
\begin{table} [ht!]
\caption{Parameters for the Heston-Cox-Ingersoll-Ross model}
\label{Table 1}
\centering
\begin{tabular}{llllllllllllll}
\hline
&\footnotesize{$\sigma $}&&\footnotesize{$\eta $}&&\footnotesize{$\rho$}&&\footnotesize{$b$}&&\footnotesize{$\alpha$}&&\footnotesize{$\beta$}\\
\hline
&\footnotesize{$0.05$}&&\footnotesize{$0.2$}&&\footnotesize{$0.8$}&&\footnotesize{$0.05$}&&\footnotesize{$0.5$}&&\footnotesize{$0.1$}\\
\hline
\end{tabular}
\vspace{.1cm}
\end{table}
\begin{table} [ht!]
\caption{The solution of HCIR PDE (\ref{(S2_18)}) at $V=0.2,R=0.2$ and time $\tau=T$ for ${k_0}={k_1}={k_2}=0$.}
\label{Table 2}
\centering
\begin{tabular}{llllllllll}
\hline
&\footnotesize{$$}&&\footnotesize{$S=120$}&&\footnotesize{$S=350$}&&\footnotesize{$S=450$}&&\footnotesize{Maximum Relative Error}\\
\hline
&\footnotesize{$M=50,J=5,K=5,N=5$}&&\footnotesize{$41.73381$}&&\footnotesize{$266.00603$}&&\footnotesize{$365.97636$}&&\footnotesize{$0.0089$}\\
&\footnotesize{$M=100,J=10,K=10,N=10$}&&\footnotesize{$42.03597$}&&\footnotesize{$266.22293$}&&\footnotesize{$366.19750$}&&\footnotesize{$0.0017$}\\
&\footnotesize{$M=200,J=20,K=20,N=10$}&&\footnotesize{$42.06774$}&&\footnotesize{$266.24294$}&&\footnotesize{$366.21802$}&&\footnotesize{$9.2290e-04$}\\
&\footnotesize{$M=800,J=40,K=40,N=15$}&&\footnotesize{$42.11332$}&&\footnotesize{$266.26775$}&&\footnotesize{$366.24379$}&&\footnotesize{$1.5959e-04$}\\
\hline
&\footnotesize{Analytical Solution \cite{G_B_1997}}&&\footnotesize{$42.10660$}&&\footnotesize{$266.30311$}&&\footnotesize{$366.29215$}\\
\hline
&\footnotesize{FD Method \cite{W_2019}}&&\footnotesize{$42.047175$}&&\footnotesize{$266.44092$}&&\footnotesize{$366.43047$}\\
\hline
\end{tabular}
\vspace{.1cm}
\end{table}
\begin{table} [ht!]
\caption{The Solution of HCIR PDE (\ref{(S2_18)}) at $V=0.2,R=0.4$ and time $\tau=T$ for ${k_0}={k_1}={k_2}=0$  }
\label{Table 22}
\centering
\begin{tabular}{llllllllll}
\hline
&\footnotesize{$$}&&\footnotesize{$S=120$}&&\footnotesize{$S=350$}&&\footnotesize{$S=450$}&&\footnotesize{Maximum Relative Error}\\
\hline
&\footnotesize{$M=50,J=5,K=5,N=5$}&&\footnotesize{$51.22891$}&&\footnotesize{$278.05398$}&&\footnotesize{$378.03803$}&&\footnotesize{$0.0106$}\\
&\footnotesize{$M=100,J=10,K=10,N=10$}&&\footnotesize{$51.55225$}&&\footnotesize{$278.23456$}&&\footnotesize{$378.22142$}&&\footnotesize{$0.0044$}\\
&\footnotesize{$M=200,J=20,K=20,N=10$}&&\footnotesize{$51.57329$}&&\footnotesize{$278.23349$}&&\footnotesize{$378.22104$}&&\footnotesize{$0.0039$}\\
&\footnotesize{$M=800,J=40,K=40,N=15$}&&\footnotesize{$51.65701$}&&\footnotesize{$278.30089$}&&\footnotesize{$378.28908$}&&\footnotesize{$0.0023$}\\
&\footnotesize{$M=1200,J=80,K=80,N=30$}&&\footnotesize{$ 51.73646$}&&\footnotesize{$278.36862$}&&\footnotesize{$378.35736$}&&\footnotesize{$7.9339e-04$}\\
\hline
&\footnotesize{Analytical Solution \cite{G_B_1997}}&&\footnotesize{$51.77754$}&&\footnotesize{$278.42940$}&&\footnotesize{$378.42374$}\\
\hline
&\footnotesize{FD Method \cite{W_2019}}&&\footnotesize{$51.55517$}&&\footnotesize{$ 278.439779$}&&\footnotesize{$378.43458$}\\
\hline
\end{tabular}
\vspace{.1cm}
\end{table}
\begin{figure} [ht!]
\begin{center}
\includegraphics[width=12cm, height=9cm]{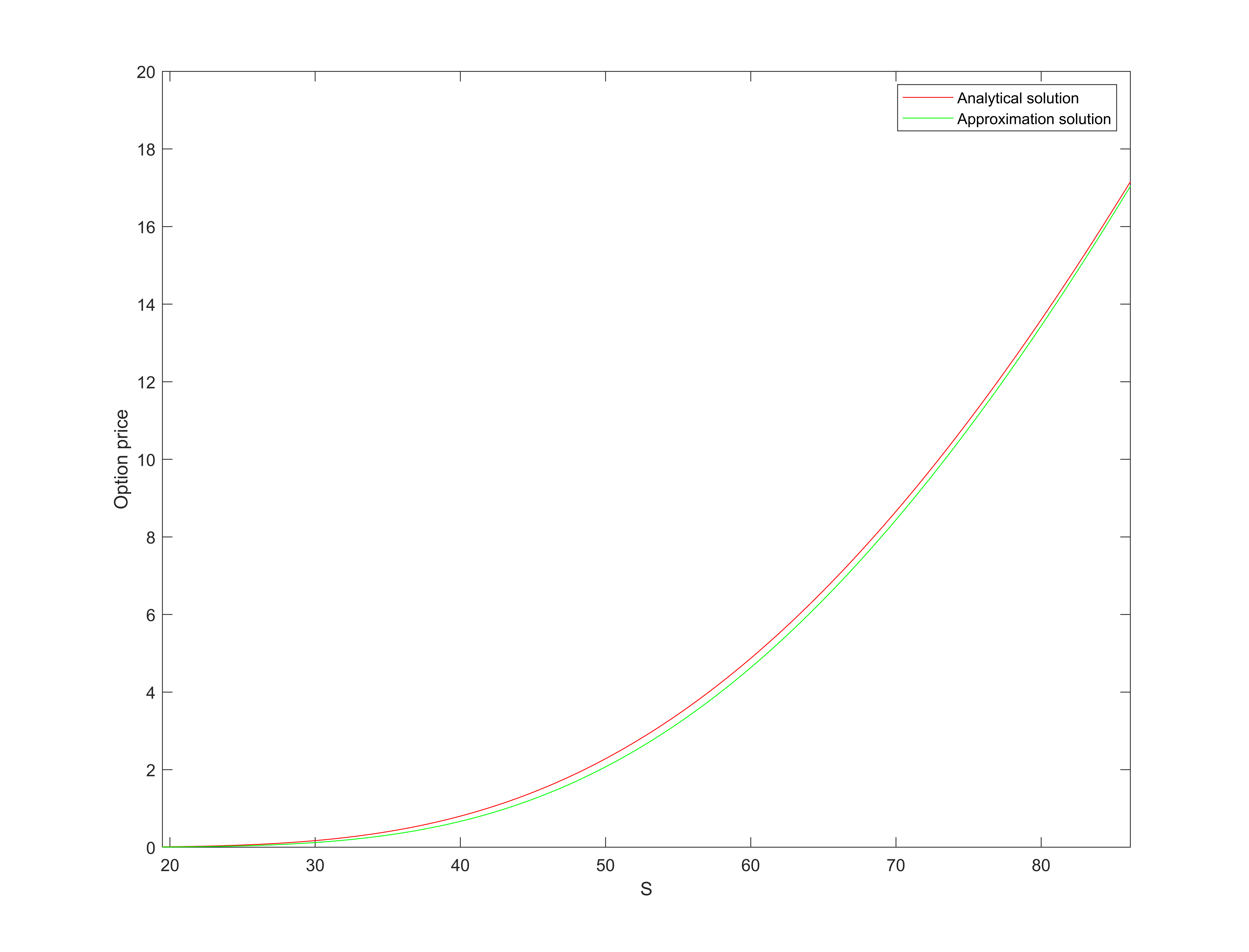}
\caption{The comparison of analytical solution and approximation solution for European call option at $\tau=T,R=0.2,V=0.2,{k_0}={k_1}={k_2}=0$.}
\label{figure1_1}
\end{center}
\end{figure}
\begin{figure} [ht!]
\begin{center}
\includegraphics[width=12cm, height=9cm]{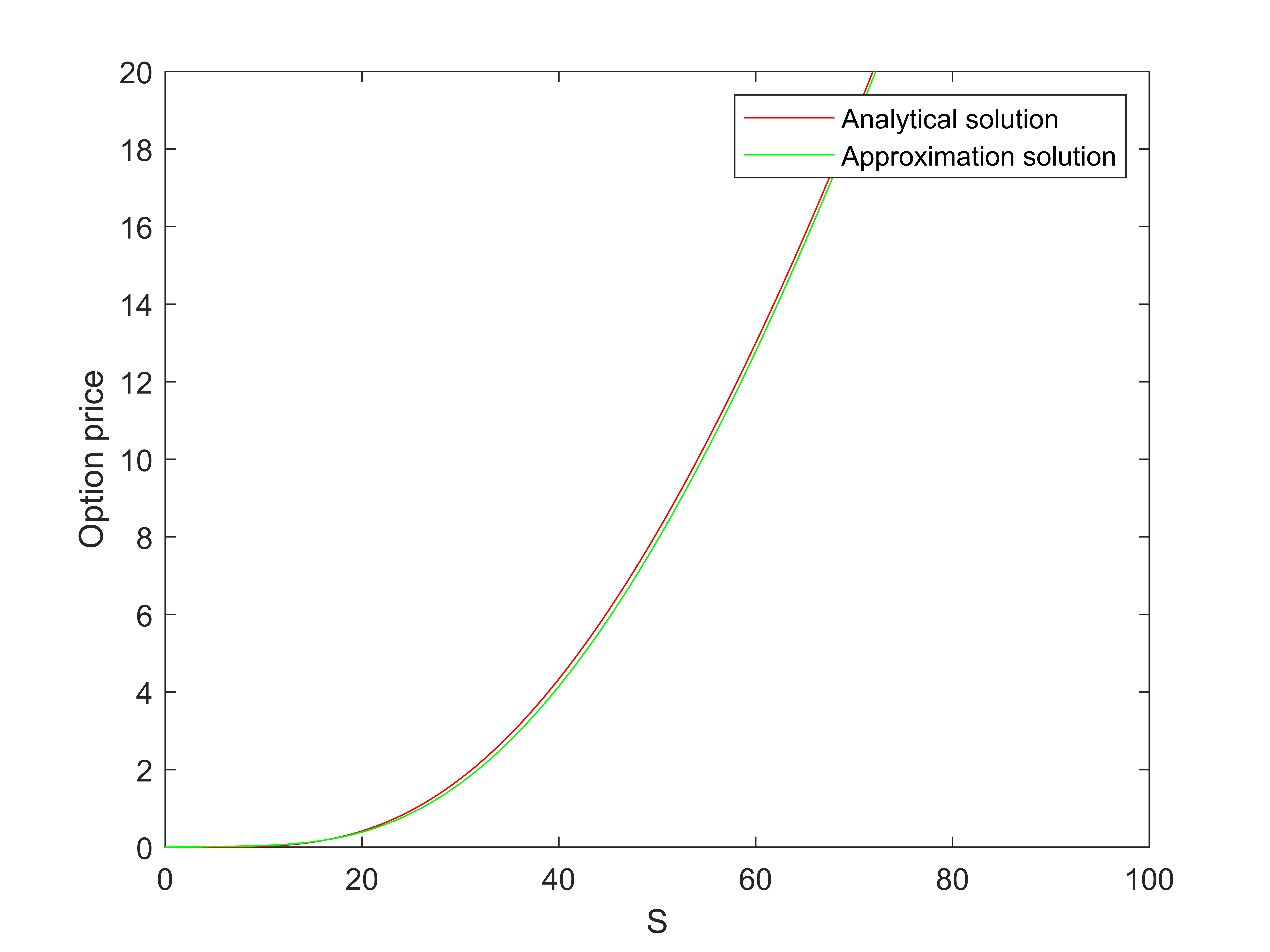}
\caption{The comparison of analytical solution and approximation solution for European call option at $\tau=T,R=0.4,V=0.4,{k_0}={k_1}={k_2}=0$.}
\label{fig2_1}
\end{center}
\end{figure}
\begin{figure} [ht!]
\begin{center}
\includegraphics[width=12cm, height=9cm]{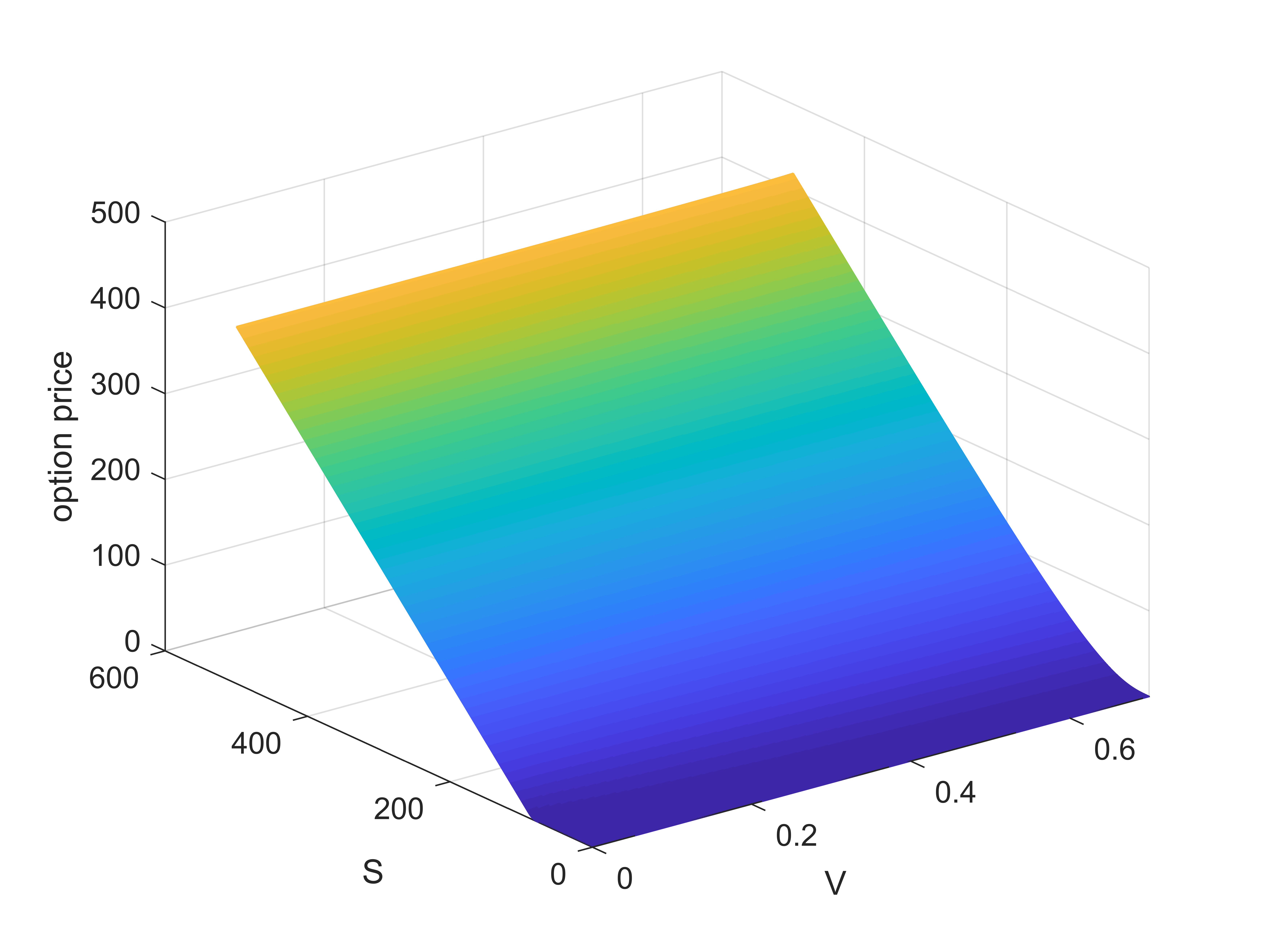}
\caption{The European call option at $\tau=T,R=0.2,{k_0}={k_1}={k_2}=0$.}
\label{figure1}
\end{center}
\end{figure}
\begin{figure} [ht!]
\begin{center}
\includegraphics[width=12cm, height=9cm]{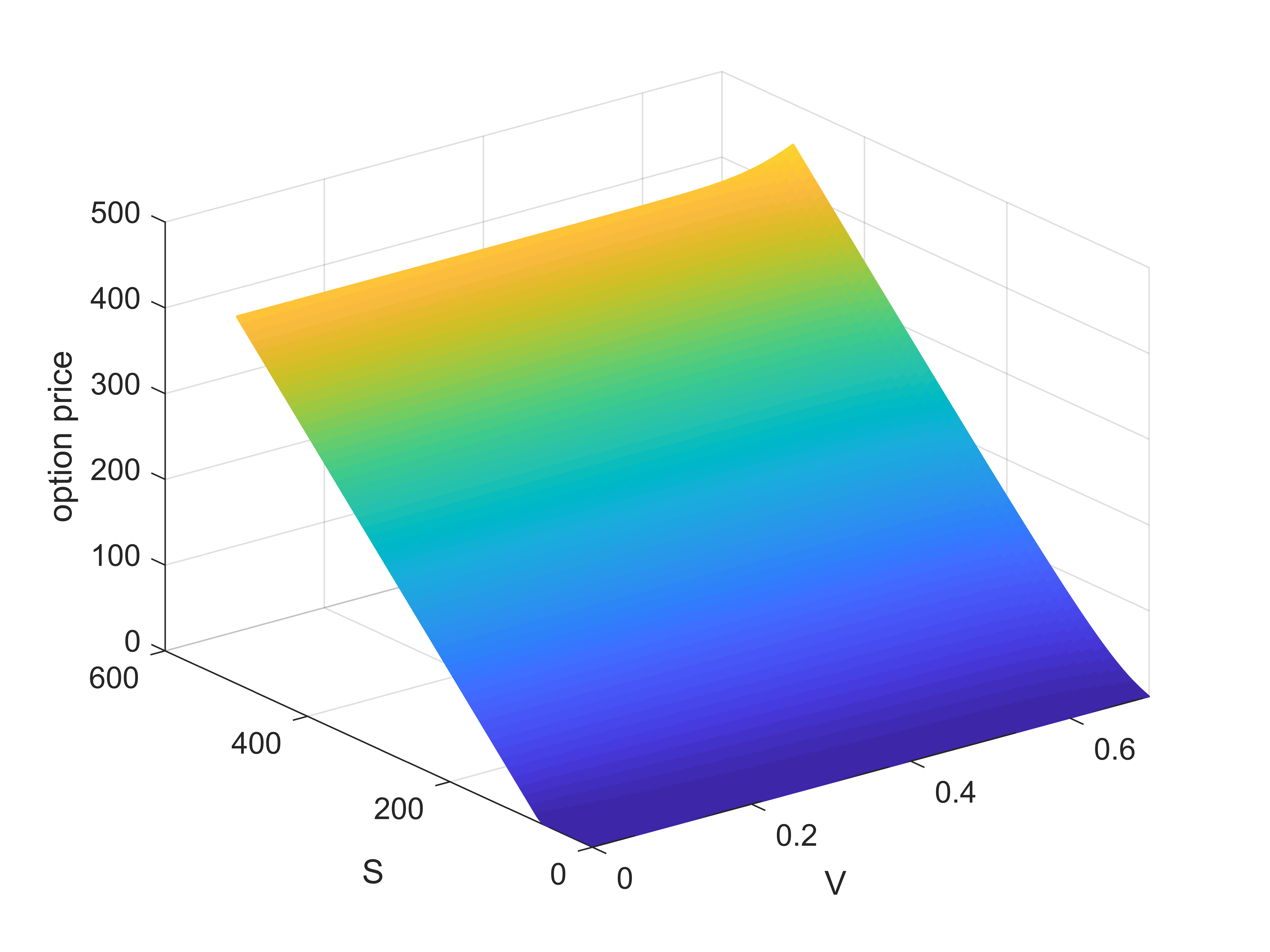}
\caption{The European call option at $\tau=T,R=0.4,{k_0}={k_1}={k_2}=0$.}
\label{fig1}
\end{center}
\end{figure}
\begin{table} [ht!]
\caption{The approximation solution of HCIR PDE (\ref{(S2_18)}) at $V=0.2,R=0.2,R=0.4$, time $\tau=T$ for ${k_0}={k_1}={k_2}=0.02$.}
\label{Table 5}
\centering
\begin{tabular}{llllllllll}
\hline
&\footnotesize{$S$}&&\footnotesize{ADI method \cite{H_2013}}&&\footnotesize{FD method \cite{W_2019}}&&\footnotesize{ADI method \cite{H_2013}}&&\footnotesize{FD method \cite{W_2019}}\\
&\footnotesize{}&&\footnotesize{$R=0.2$}&&\footnotesize{$ R=0.2$}&&\footnotesize{$R=0.4$}&&\footnotesize{$ R=0.4$}\\
\hline
&\footnotesize{$120$}&&\footnotesize{$41.61569$}&&\footnotesize{$39.67547$}&&\footnotesize{$51.28916$}&&\footnotesize{$51.74822$}\\
&\footnotesize{$350$}&&\footnotesize{$266.07274$}&&\footnotesize{$264.33191$}&&\footnotesize{$278.12679$}&&\footnotesize{$279.66139$}\\
&\footnotesize{$450$}&&\footnotesize{$366.05583$}&&\footnotesize{$364.32657$}&&\footnotesize{$378.11851$}&&\footnotesize{$379.65636$}\\
\hline
\end{tabular}
\vspace{.1cm}
\end{table}
\begin{figure} [ht!]
\begin{center}
\includegraphics[width=12cm, height=9cm]{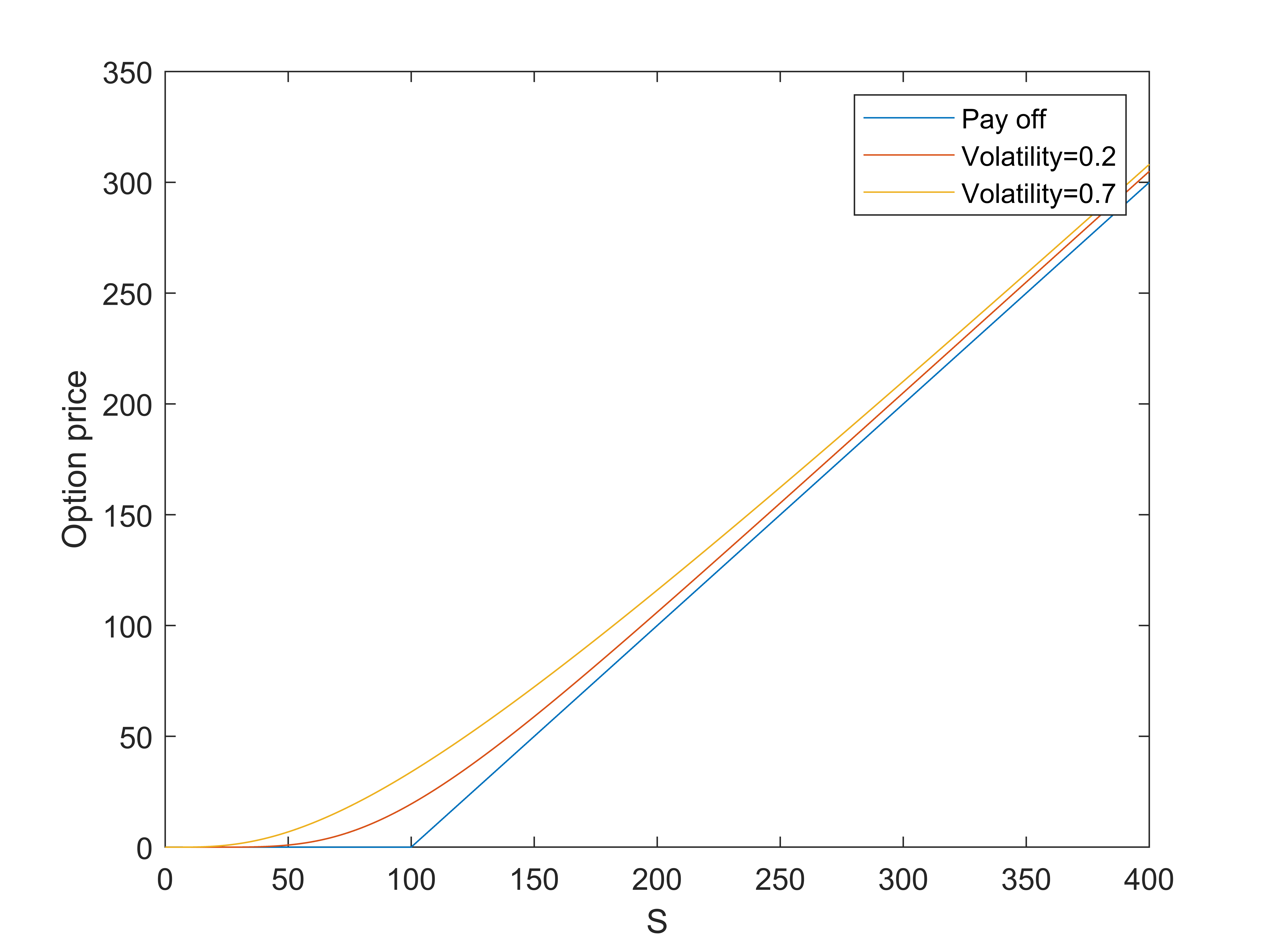}
\caption{The European call option at $\tau=T,R=0.2,V=0.2,V=0.7$ and ${k_0}={k_1}={k_2}=0.02$.}
\label{figure2}
\end{center}
\end{figure}
\begin{figure} [ht!]
\begin{center}
\includegraphics[width=12cm, height=9cm]{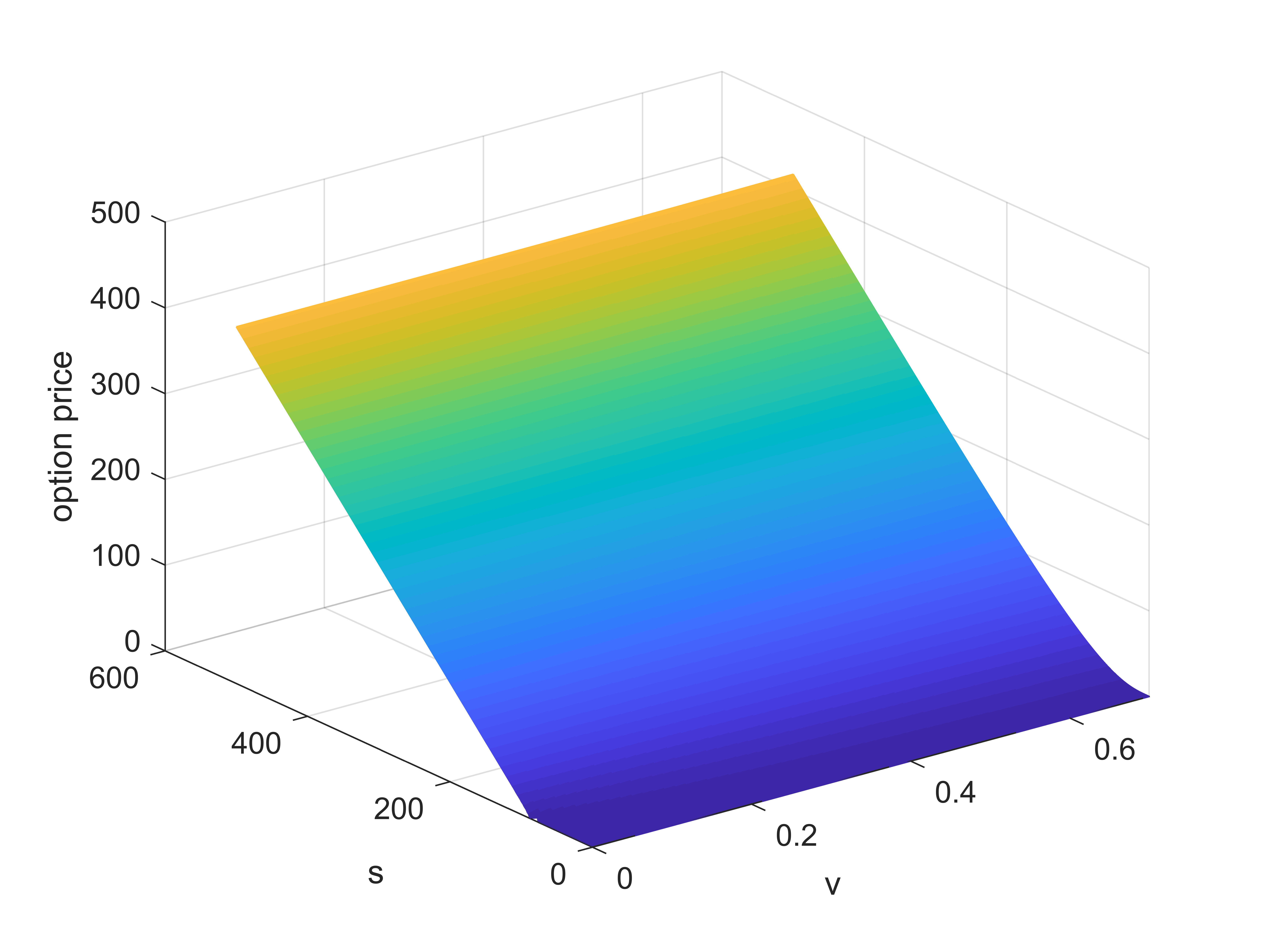}
\caption{The European call option at $\tau=T,R=0.2$ and ${k_0}={k_1}={k_2}=0.02$.}
\label{figure22}
\end{center}
\end{figure}
\begin{figure} [ht!]
\begin{center}
\includegraphics[width=12cm, height=9cm]{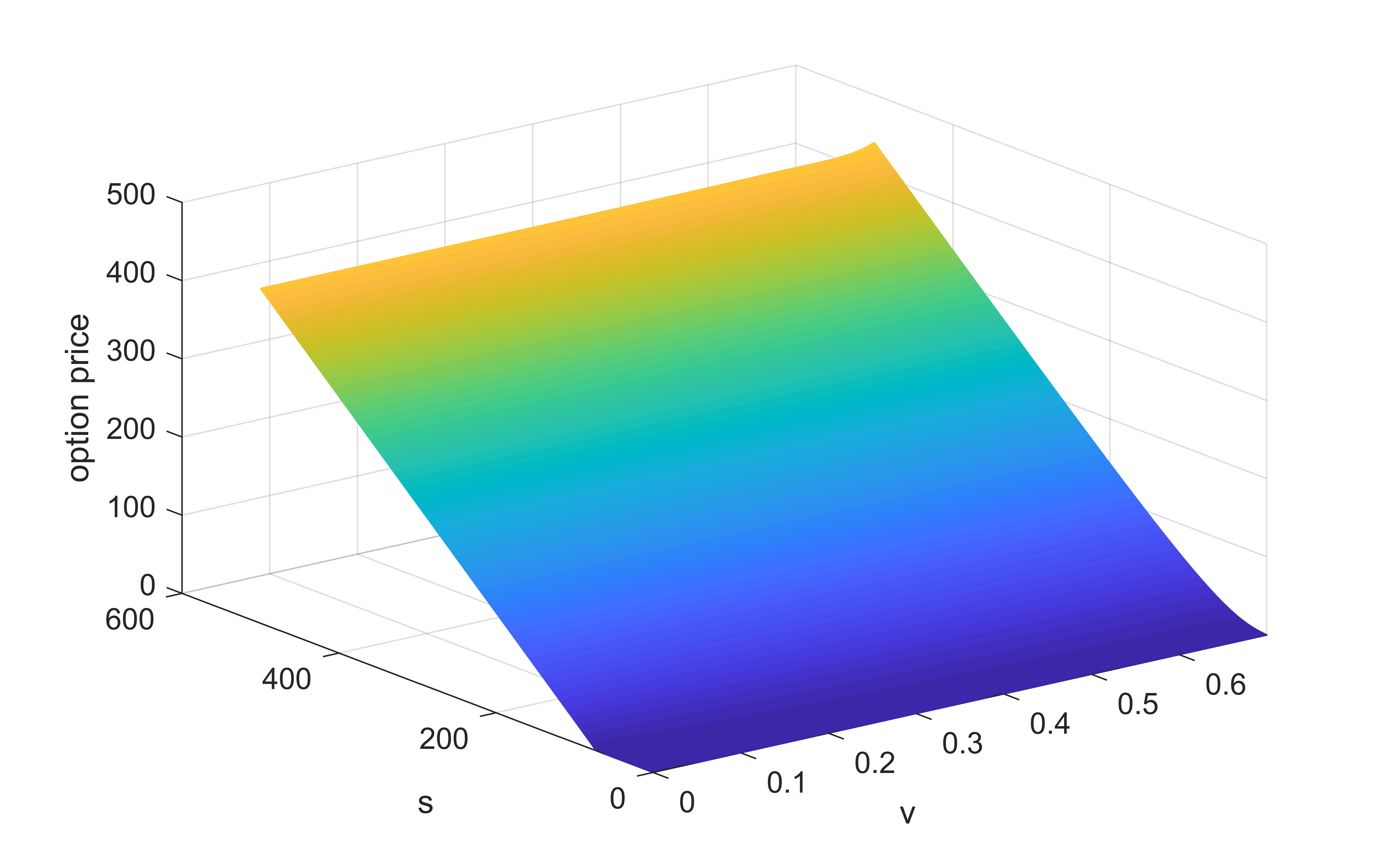}
\caption{The European call option at $\tau=T,R=0.4$ and ${k_0}={k_1}={k_2}=0.02$.}
\label{figure222}
\end{center}
\end{figure}
\begin{figure} [ht!]
\begin{center}
\includegraphics[width=12cm, height=9cm]{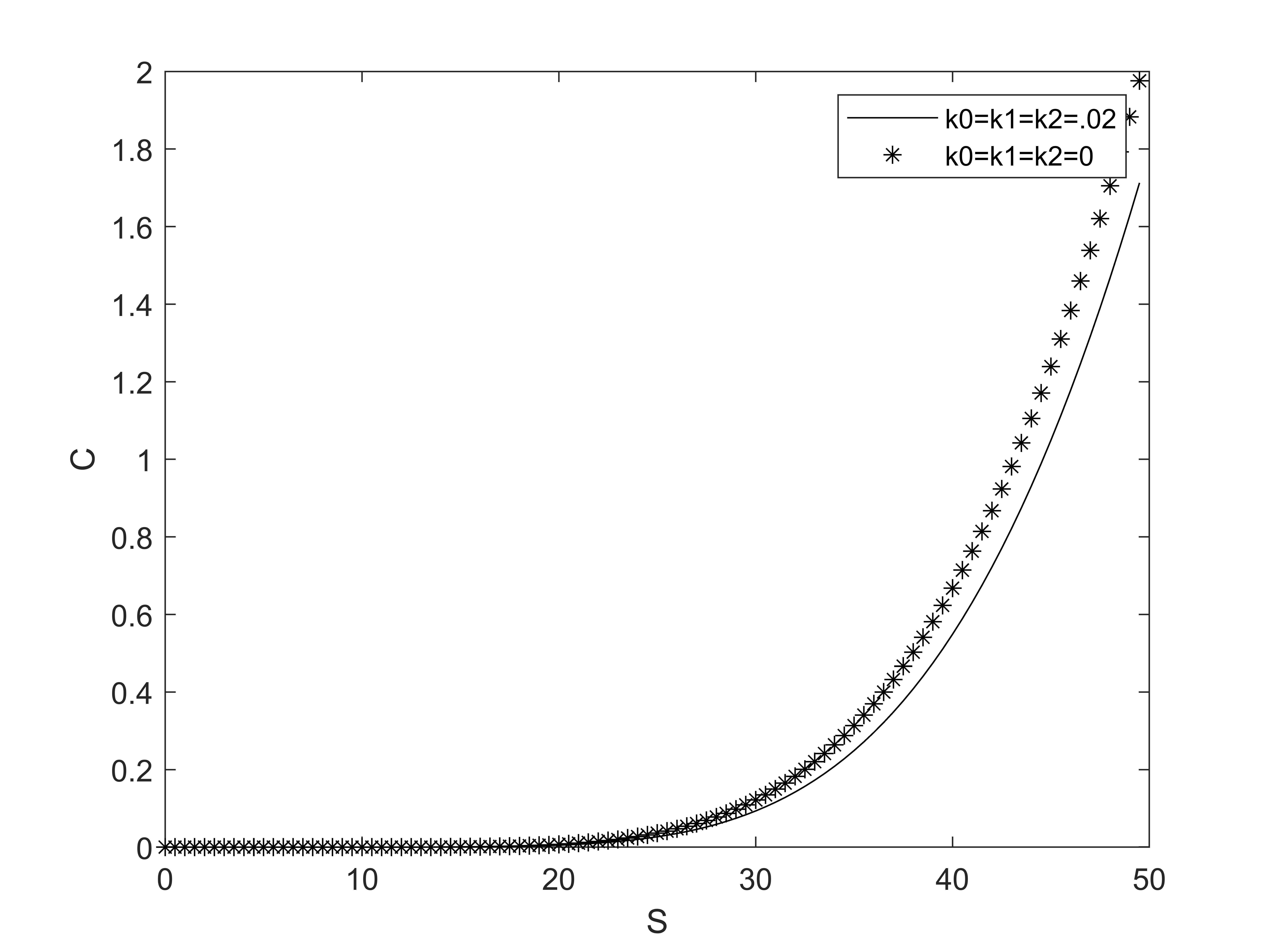}
\caption{The European call option at $\tau=T$, $R=0.2$, $V=0.2$ and ${k_0}={k_1}={k_2}=0.02$ and ${k_0}={k_1}={k_2}=0$.}
\label{figure3}
\end{center}
\end{figure}
\\ \\ \\ \\ \\ \\ \\  \\ \\ \\ \\ \\ \\ \\ \\ \\ \\
\section{Conclusion}
\label{section6}
European call contract option prices and transaction costs were gathered in this study using the HCIR model, which is a combination of the Heston \cite{Heston_1993} and Cox-Ingersoll-Ross \cite{C_I_R_1985} models. Using the Douglas scheme \cite{H_2013,D_1962}, which is a type of alternating direction implicit method, we solved the HCIR PDE problem and conducted a numerical investigation of convergence. The stability of the ADI strategy was acted in \cite{H_K_2013} utilizing Von Neumann. We compared the price of a European call option with and without transaction costs, demonstrating that the European call option price with transaction costs, has lower values than the European call option price without transaction costs in most points of the grid.\\ \\ \\
\textbf{Declarations}\\\\
   The authors declare that they have no conflict of interest.

\end{quotation}

\begin{thebibliography}{99}
\bibitem{B_1973} Black, F., Scholes, M.: The Pricing of Options and Corporate Liabilities. J. Political Econ. 81, 637-654 (1973)
\bibitem{H_W_1987} Hull, J., White, A.: The Pricing of Options on Assets with Stochastic Volatilities. J. Financ. 42, 281-300 (1987) 
\bibitem{S_Z_1999} Sch$\ddot o$bel,  R., Zhu, J.: Stochastic Volatility With an Ornstein-Uhlenbeck Process: An Extension. Rev. Financ. 3, 23-46 (1999)
\bibitem{S_S_1991} Stein, E.M., Stein, J.C.: Stock Price Distributions with Stochastic Volatility: An Analytic Approach. Rev. Financ. Stud. 4, 727-752 (1991) 
\bibitem{Heston_1993} Heston, S.L.: A Closed-Form Solution for Options with Stochastic Volatility with Applications to Bond and Currency Options. Rev. Financ. Stud. 6, 327-343 (1993) 
\bibitem{G_O_W_2011} Grzelak, L.A., Oosterlee, C.W., Van Weeren, S.: The affine Heston model with correlated Gaussian interest rates for pricing hybrid derivatives. Quant. Finance. 11, 1647-1663 (2011)
\bibitem{G_O_W_2012} Grzelak, L.A., Oosterlee, C.W., Van Weeren, S.: Extension of stochastic volatility equity models with the Hull-White interest rate process. Quant. Financ. 12, 89-105 (2012) 
\bibitem{G-G_O_2013} Guo, S., Grzelak, L.A., Oosterlee, C.W.: Analysis of an affine version of the Heston-Hull-White option pricing partial differential equation. Appl. Numer. Math. 72,  143-159 (2013)
\bibitem{H_L_P_S_2000} Van Haastrecht, A., Lord, R., Pelsser, A., Schrager, D.: Pricing long-dated insurance contracts with stochastic interest rates and stochastic volatility. Insur. Math. Econ. 45, 436-448 (2009) 
\bibitem{C_I_R_1985} Cox, J.C., Ingersoll, J.E., Ross, S.A.: An intertemporal general equilibrium model of asset prices. Econometrica. 53, 385-407 (1985) 
\bibitem{H_W_1990} Hull, J., White, A.: Pricing Interest-Rate-Derivative Securities. Rev. Financ. Stud. 3, 573-592 (1990)
\bibitem{V_1977} Vasicek, O.: An equilibrium characterization of the term structure. J. Financ. Econ. 5, 177-188 (1977)
\bibitem{S_O_2002} Sippel, J., Ohkoshi, S.: All power to PRDC notes, Risk, 15, 531-533 (2002) 
\bibitem{G_O_2011} Grzelak, L.A., Oosterlee, C.W.: On the Heston Model with Stochastic Interest Rates. J. Financ. Math. 2, 255-286 (2011) 
\bibitem{D_2006} Duffy, Daniel J.: Finite difference methods in financial engineering. John Wiley and Sons Ltd, Chichester (2006)
\bibitem{W_2019} Wang, B.: Option pricing under the Heston-CIR model with stochastic interest rates and transaction costs. PhD diss., Auckland University of Technology, New Zealand (2019)
\bibitem{C_W_Z_2021} Cao, J., Wang, B., Zhang, W.: Valuation of European options with stochastic interest rates and transaction costs. Int. J. Comput. Math. 99(2), 227-239 (2022)
\bibitem{H_2013} Haentjens, T.: Efficient and stable numerical solution of the Heston-Cox-Ingersoll-Ross partial differential equation by alternating direction implicit finite difference schemes, Int. J. Comput. Math. 90, 2409-2430 (2013)
\bibitem{Z_2002} Zapart, Christopher.: Stochastic volatility options pricing with wavelets and artificial neural networks. Quantitative Finance 2(6), 487 (2002)
\bibitem{N_Y_2011} Niederreiter, H.: Monte Carlo and Quasi-Monte Carlo Methods 1998: Proceedings of a Conference Held at the Claremont Graduate University, Claremont, California, USA, JU. Springer-Verlag (2000)
\bibitem{Douglas} Douglas, J., and Rachford, H.H.: on the numerical solution of heat conduction equations in two and three dimensions. Trans. Am. Math. Assoc. 82, 421-439 (1955)
\bibitem{D_1955} Douglas, J.: On the numerical integration of ${u_{xx}} + {u_{yy}} = {u_t}$ by implicit methods. J. Soc. Ind. Appl. Math. 3, 42-65 (1955) 
\bibitem{Peaceman} Peaceman, D.W.: "Differential equations for flow in reservoirs." Fundamentals of numerical reservoir simulation. Amsterdam, Elsevier 6, 1-34 (1977) 
\bibitem{P_R_1955} Peaceman, D.W., and Rachford, H.H.: The numerical solution of parabolic and elliptic differential equations, J. Soc. Ind. Appl. Math. 3, 28-41 (1955) 
\bibitem{D_R_1956} Douglas,  J., and Rachford, H.H.: On the numerical solution of heat conduction problems in two and three space variables. Trans. Amer. Math. Soc. 82, 421-439 (1956)
\bibitem{Y_1971} Yanenko, N.N.: The Method of Fractional Steps, The Solution of Problems of Mathematical Physics in Several Variables. Heidelberg, Springer Berlin (1971).
\bibitem{D_1962} Douglas, J.: Alternating Direction Methods for Three Space Variables, Numerische Mathematik 4, 41-63 (1962)  
\bibitem{S_2014} Shidfar, A., Paryab, K., Yazdanian, A.R., Pirvu, T.A.: Numerical analysis for Spread option pricing model of markets with finite liquidity: first-order feedback model. Int. J. Comput. Math. 91, 2603-2620 (2014)
\bibitem{Y_2014} Yazdanian, A.R. and Pirvu, T.A.: Numerical analysis for Spread option pricing model in illiquid underlying asset market: full feedback model. Quant. Finance. (2014) https://doi.org/10.48550/arXiv.1406.1149 
\bibitem{D_F_2012} D$\ddot u$ring, B., Fourni$e'$, M.: High-order compact finite difference scheme for option pricing in stochastic volatility models. J. Comput. Appl. Math. 236, 4462-4473 (2012)
\bibitem{S_N_N_2019} Safaei, M., Neisy, A., Nematollahi, N.: Generalized Componentwise Splitting Scheme For Option Pricing Under The Heston-Cox-Ingersoll-Ross Model. Journal of Statistical Theory and Applications 18, 425-438 (2019)
\bibitem{J_S_M_W_2022} Wang, J., Wen, S., Yang, M., Shao, W.: Practical finite difference method for solving multi-dimensional Black-Scholes model in fractal market, Chaos, Solitons and Fractals 157, 111895 (2022)
\bibitem{Leland_1985} Leland, H.: Option pricing and replication with transactions costs. The journal of finance 40(5), 1283-1301 (1985)
\bibitem{H_N_1989} Hodges, S., and Neuberger, A.: Optimal replication of contingent claims under transaction costs. Review Futures Market 8, 222-239 (1989) 
\bibitem{G_S_1996} Grannan, E., and Swindle, G.: Minimizing transaction costs of option hedging strategies. Mathematical finance 6(4), 341-364 (1996)
\bibitem{Z_Z_2007} Zhao, Y., and Ziemba, W.: Comments on and corrigendum to “Hedging errors with Leland’s option model in the presence of transaction costs”, Financ. Res. Lett. 4, 196-199 (2007) 
\bibitem{Zh_Z_2007} Zhao, Y., and Ziemba, W.: Hedging errors with Leland’s option model in the presence of transaction costs. Finance. Res. Lett. 4, 49-58 (2007) 
\bibitem{Se_2014} SenGupta, I.: Option pricing with transaction costs and stochastic interest rate. Appl. Math. Financ. 21, 399-416 (2014) 
\bibitem{M_2015} Mariani, M.C., SenGupta, I., and Sewell, G.: Numerical methods applied to option pricing models with transaction costs and stochastic volatility. Quant. Financ. 15, 1417-1424 (2015)
\bibitem{N_P_2015} Nguyen, T.H., and Pergamenschchikov, S.: Approximate hedging with proportional transaction costs in stochastic volatility models with jumps. Mathematical Finance (2020). https://doi.org/10.48550/arXiv.1505.02627
\bibitem{N_P_2017} Nguyen, T.H., and Pergamenshchikov, S.: Approximate hedging problem with transaction costs in stochastic volatility markets, Math. Financ. 27, 832-865 (2017)
\bibitem{G_B_1997} Bakshi, G., Cao, C., and Chen, Z.: Empirical performance of alternative option pricing models. Journal of Finance 52, 2003-2049 (1997) 
\bibitem{C_J_L_2016} Cao, J., Lian, G., and Roslan, T.R.N.: Pricing variance swaps under stochastic volatility and stochastic interest rate. Appl. Math. Comput. 277, 72-81 (2016)
\bibitem{H_K_2013}  in't Hout, K.J., and Mishra, C.: Stability of ADI schemes for multidimensional diffusion equations with mixed derivative terms, Appl. Numer. Math. 74, 83-94 (2013) 
\bibitem{SH_2004} Shreve, S.E.: Stochastic calculus for finance II: Continuous-time models. Springer, New York (2004)
\end{thebibliography}
\end{document}